# SOME REMARKS ON POINTS OF LEBESGUE DENSITY AND DENSITY-DEGREE FUNCTIONS

SILVANO DELLADIO

ABSTRACT. Some properties of $m$-density points and density-degree functions are studied. Moreover the following main results are provided:

- Let $\lambda$ be a continuous differential form of degree $h$ in $\mathbb{R}^n$ (with $h \geq 0$) having the following property: There exists a continuous differential form $\Delta$ of degree $h+1$ in $\mathbb{R}^n$ such that
$$\int_{\mathbb{R}^n} \Delta \wedge \omega = \int_{\mathbb{R}^n} \lambda \wedge d\omega,$$
for every $C_c^\infty$ differential form $\omega$ of degree $n-h-1$ in $\mathbb{R}^n$. Moreover let $\mu$ be a $C^1$ differential form of degree $h+1$ in $\mathbb{R}^n$ and set $E := \{y \in \mathbb{R}^n \,|\, \Delta(y) = \mu(y)\}$. Then $d\mu(x) = 0$ whenever $x$ is a $(n+1)$-density point of $E$.

- Let $f : \mathbb{R}^n \to \overline{\mathbb{R}}$ be a measurable function such that $f(x) \in \{0\} \cup [n, +\infty]$ for a.e. $x \in \mathbb{R}^n$. Then there exists a countable family $\{F_k\}_{k=1}^\infty$ of closed subsets of $\mathbb{R}^n$ such that the corresponding sequence of density-degree functions $\{d_{F_k}\}_{k=1}^\infty$ converges almost everywhere to $f$.

## 1. INTRODUCTION

Let $m \in [n, +\infty)$ and $E \subset \mathbb{R}^n$. Then $x \in \mathbb{R}^n$ is said to be an $m$-density point of $E$ if $\mathcal{L}^n(B(x,r) \setminus E) = o(r^m)$ as $r \to 0^+$, where $\mathcal{L}^n$ denotes the Lebesgue outer measure in $\mathbb{R}^n$ and $B(x,r)$ is the open ball of radius $r$ centered at $x$. The set of $m$-density points of $E$ is denoted by $E^{(m)}$ (cf. [2, 3]). For $x \in E^{(n)}$, the supremum of the numbers $m$ such that $x \in E^{(m)}$ is called the density-degree of $E$ at $x$ and is denoted by $d_E(x)$, while $d_E(x)$ is assumed to be zero when $x \notin E^{(n)}$ (cf. [7]). Thanks to these definitions, points of Lebesgue density cease to be indistinguishable from each other and are instead characterised by their own density-degree.

There are reasons for considering points with a high density-degree, albeit finite, as interior points. For example, the following fact holds: If $x \in \{\operatorname{grad} f = F\}^{(n+1)}$, where $f \in C^1(\mathbb{R}^n)$ and $F \in C^1(\mathbb{R}^n, \mathbb{R}^n)$, then the Jacobian matrix $DF(x)$ is symmetric. We thus discover that the geometry of $\{\operatorname{grad} f = F\}$ is characterised by a very low density at points where $DF$ is non-symmetric (cf. [2, Theorem 2.1]; subsequent extensions to the context of PDE can be found e.g. in [10, Theorem 3.6] and [11, Theorem 3.2]).

A similar application can be given in the context of the Frobenius theorem on distributions. More precisely let $M$ and $D$ be, respectively, a $n$-dimensional $C^1$ submanifold of $\mathbb{R}^{n+k}$ and a $C^1$ distribution of rank $n$ on $\mathbb{R}^{n+k}$. Then $D$ is involutive at each $(n+1)$-density point of the tangency set of $M$ with respect to $D$. Hence the tangency must be low in density at all points where the distribution $D$ is non-involutive (cf. [8, Theorem 1.1]).

Another fact worth mentioning in this introduction is the following: Except for a subset of null measure, the points of a set of locally finite perimeter are $m$-density points, with $m := n + 1 + \frac{1}{n-1}$ (cf. [2, Lemma 4.1]; for further applications in this context, see [5, 6]). Using $m$-density one can also define the notion of $m$-approximate continuity, which for $m = n$ reduces to the well known approximate continuity (cf. [13, Section 2.9.12] and [12, Section 1.7.2]). Properties of $m$-approximate continuity holds for Sobolev functions and for functions of bounded variation (cf. [9]).

Finally, let us briefly describe the two main results of this note. The first one, proved in Section 3 below, is the following generalization of [2, Theorem 2.1].

**Theorem.** *Let $\lambda$ be a continuous differential form of degree $h$ in $\mathbb{R}^n$ (with $h \geq 0$) having the following property: There exists a continuous differential form $\Delta$ of degree $h+1$ in $\mathbb{R}^n$ such that*

$$\int_{\mathbb{R}^n} \Delta \wedge \omega = \int_{\mathbb{R}^n} \lambda \wedge d\omega,$$

*for every $C_c^\infty$ differential form $\omega$ of degree $n - h - 1$ in $\mathbb{R}^n$. Moreover let $\mu$ be a $C^1$ differential form of degree $h + 1$ in $\mathbb{R}^n$ and set $E := \{y \in \mathbb{R}^n \,|\, \Delta(y) = \mu(y)\}$. Then $(d\mu)|_{E^{(n+1)}} \equiv 0$.*

The second result is provided in Section 5. It originates from the following question: Is it true that for every (measurable) function $f : \mathbb{R}^n \to \{0\} \cup [n, +\infty]$ there exists $E \subset \mathbb{R}^n$ such that $d_E = f$ almost everywhere? Some hasty considerations may mislead us into thinking that the answer is yes, but on deeper reflection, it is not difficult to conclude that the correct answer is no (cf. Example 5.1 below). However, somewhat surprisingly, the following approximation property holds true.

**Theorem.** *Let $f : \mathbb{R}^n \to \overline{\mathbb{R}}$ be a measurable function such that $f(x) \in \{0\} \cup [n, +\infty]$ for a.e. $x \in \mathbb{R}^n$. Then there exists a countable family $\{F_k\}_{k=1}^\infty$ of closed subsets of $\mathbb{R}^n$ such that $\lim_{k \to \infty} d_{F_k}(x) = f(x)$ for a.e. $x \in \mathbb{R}^n$.*

## 2. Preliminaries

**2.1. General notation.** The coordinates of $\mathbb{R}^n$ are denoted by $(x_1, \ldots, x_n)$ and we set $D_i := \partial/\partial x_i$. If $k$ is any positive integer not exceeding $n$, then $I(n,k)$ is the family of integer multi-indices $\alpha = (\alpha_1, \ldots, \alpha_k)$ such that $1 \leq \alpha_1 < \cdots < \alpha_k \leq n$. If $\alpha \in I(n,k)$, then we denote by $\bar{\alpha}$ the member of $I(n, n-k)$ which complements $\alpha$ in $\{1, 2, \ldots, n\}$ in the natural increasing order (e.g., if $\alpha = (2,3,5) \in I(7,3)$, then $\bar\alpha = (1,4,6,7) \in I(7,4)$). The open ball of radius $r$ centered at $x \in \mathbb{R}^n$ is denoted by $B(x,r)$. Sometimes, for simplicity, the sphere $B(0,r)$ will be denoted by $B_r$. Moreover $\chi_E$ is the characteristic function of $E$. The equivalence relation of functions and the equivalence relation of subsets of $\mathbb{R}^n$, with respect to the Lebesgue outer measure $\mathcal{L}^n$, are both denoted by $\sim$. If $E, F \subset \mathbb{R}^n$ and $\mathcal{L}^n(E \setminus F) = 0$, then we write $E \widetilde\subset F$ (so that $E \sim F$ if and only if $E \widetilde\subset F$ and $F \widetilde\subset E$). Observe that $E \widetilde\subset F$ if and only if $F^c \widetilde\subset E^c$.

2.2. **Covectors and differential forms in $\mathbb{R}^n$.** If $k$ is a positive integer not exceeding $n$, then a $k$-covector (of $\mathbb{R}^n$) is a $k$-linear alternating map from $(\mathbb{R}^n)^k$ to $\mathbb{R}$. Let $\bigwedge^k(\mathbb{R}^n)$ denote the set of all $k$-covectors. In particular $\bigwedge^1(\mathbb{R}^n)$ is the dual space of $\mathbb{R}^n$ and we will denote the standard dual basis by $dx_1, \ldots, dx_n$. The set $\bigwedge^k(\mathbb{R}^n)$ is a vector space of dimension $\binom{n}{k}$ with the standard basis

(2.1) $$\{dx_\alpha := dx_{\alpha_1} \wedge \cdots \wedge dx_{\alpha_k} \,|\, \alpha \in I(n,k)\}$$

where $\wedge$ denotes the wedge product. Recall that $\bigwedge^k(\mathbb{R}^n)$ is equipped with the following inner product naturally induced from $\mathbb{R}^n$ (making (2.1) an orthonormal basis):

$$\langle \xi, \eta \rangle := \sum_{\alpha \in I(n,k)} \xi_\alpha \eta_\alpha$$

with

$$\xi = \sum_{\alpha \in I(n,k)} \xi_\alpha \, dx_\alpha, \quad \eta = \sum_{\alpha \in I(n,k)} \eta_\alpha \, dx_\alpha.$$

We observe that the following property holds (we set for simplicity $dx := dx_1 \wedge \cdots \wedge dx_n$):

(2.2) If $\xi \in \bigwedge^k(\mathbb{R}^n)$ and $\langle \xi \wedge \eta, dx \rangle = 0$ for all $\eta \in \bigwedge^{n-k}(\mathbb{R}^n)$, then $\xi = 0$.

Indeed, if $\xi = \sum_{\alpha \in I(n,k)} \xi_\alpha \, dx_\alpha$, then for all $\beta \in I(n,k)$ we have $0 = \langle \xi \wedge dx_{\bar\beta}, dx \rangle = \xi_\beta \langle dx_\beta \wedge dx_{\bar\beta}, dx \rangle$, hence $\xi_\beta = 0$.

A $C^H$ differential form of degree $k$ (on $\mathbb{R}^n$) is a $k$-covector field

$$x \in \mathbb{R}^n \mapsto \sum_{\alpha \in I(n,k)} f_\alpha(x) \, dx_\alpha,$$

with $\{f_\alpha \,|\, \alpha \in I(n,k)\} \subset C^H(\mathbb{R}^n)$. A $C^H$ differential form of degree 0 is simply a function in $C^H(\mathbb{R}^n)$. We recall that the addition and the exterior product of covectors naturally induce the addition and the exterior product of differential forms. Moreover an exterior derivative operator $d$ is uniquely defined on $C^1$ differential forms and it holds that

- On $C^1$ differential forms of degree 0, the operator $d$ agrees with the ordinary differential on $C^1(\mathbb{R})$;
- $d(\lambda + \mu) = d\lambda + d\mu$, for all $C^1$ differential forms $\lambda$ and $\mu$ of degree $k$;
- $d(\lambda \wedge \mu) = d\lambda \wedge \mu + (-1)^k \lambda \wedge d\mu$, for all $C^1$ differential forms $\lambda$ and $\mu$, if $\lambda$ has degree $k$;
- $d(d\lambda) = 0$, for all $C^2$ differential forms $\lambda$.

We also recall that, given a continuous differential form $\omega$ of degree $n$ with compact support and a measurable set $E \subset \mathbb{R}^n$, the integral of $\omega$ on $E$ is defined as follows:

$$\int_E \omega := \int_E \langle \omega(y), dx_1 \wedge \cdots \wedge dx_n \rangle \, d\mathcal{L}^n(y).$$

For a comprehensive treatment of $k$-covectors and differential forms, we refer the reader to the numerous books on differential geometry and geometric measure theory that deal with this subject, e.g., [16] and [15].

2.3. **Points of density.** We recall the definition of $m$-density point (compare [2, 3, 4]).

**Definition 2.1.** *Let $m \in [n, +\infty)$ and $E \subset \mathbb{R}^n$. Then $x \in \mathbb{R}^n$ is said to be a "$m$-density point of $E$" if*

$$\lim_{r \to 0^+} \frac{\mathcal{L}^n(B(x,r) \setminus E)}{r^m} = 0.$$

*The set of $m$-density points of $E$ is denoted by $E^{(m)}$.*

**Remark 2.1.** *The following simple facts occur:*

(1) *Every interior point of $E \subset \mathbb{R}^n$ is an $m$-density point of $E$, for all $m \in [n, +\infty)$. Thus, whenever $E$ is open, one has $E \subset E^{(m)}$ for all $m \in [n, +\infty)$.*
(2) *If $n \leq m_1 \leq m_2 < +\infty$ and $E \subset \mathbb{R}^n$, then $E^{(m_2)} \subset E^{(m_1)}$. In particular, one has $E^{(m)} \subset E^{(n)}$ for all $m \in [n, +\infty)$.*
(3) *Let $\{E_j\}_{j \in J}$ be any family of subsets of $\mathbb{R}^n$ and $m \in [n, +\infty)$.*
    – *One has*

$$\left( \bigcap_{j \in J} E_j \right)^{(m)} \subset \bigcap_{j \in J} E_j^{(m)}, \quad \left( \bigcup_{j \in J} E_j \right)^{(m)} \supset \bigcup_{j \in J} E_j^{(m)}$$

    – *Let $J$ be finite. Then*

(2.3)
$$\left( \bigcap_{j \in J} E_j \right)^{(m)} = \bigcap_{j \in J} E_j^{(m)},$$

    *while the identity*

$$\left( \bigcup_{j \in J} E_j \right)^{(m)} = \bigcup_{j \in J} E_j^{(m)}$$

    *fails to be true in general, e.g., $E_1 = (-1, 0)$ and $E_2 = (0, 1)$ (with $n = 1$ and $J = \{1, 2\}$).*
    – *If $J$ is countably infinite, then (2.3) can fail to be true, e.g., $J = \{j = 1, 2, \ldots\}$ and $E_j := B(0, 1/j)$.*
(4) *If $E \subset \mathbb{R}^n$ is measurable, then $(E^{(n)})^{(m)} = E^{(m)}$ for all $m \in [n, +\infty)$. In particular $(E^{(n)})^{(n)} = E^{(n)}$.*

**Proposition 2.1** ([6], Proposition 3.1). *For all $E \subset \mathbb{R}^n$, the set $E^{(m)}$ is measurable.*

**Theorem 2.1** ([7], Corollary 4.1). *If $E$ is a measurable subset of $\mathbb{R}^n$ and $m \in (n, +\infty)$, then*
$$E^{(m)} \sim \left\{ x \in \mathbb{R}^n \,\middle|\, \limsup_{r \to 0^+} \frac{\mathcal{L}^n(B(x,r) \setminus E)}{r^m} < +\infty \right\}.$$

The Lebesgue density theorem states that if $E$ is a measurable subset of $\mathbb{R}^n$, then almost every $x \in E$ is a $n$-density point of $E$. A remarkable family of sets that turn out to be strictly more dense than generic measurable sets is that of finite perimeter sets. We recall that the perimeter of a measurable set $E \subset \mathbb{R}^n$, denoted by $P(E)$, is the variation of $\chi_E$, that is
$$P(E) := \sup \left\{ \int_E \operatorname{div} \varphi \, d\mathcal{L}^n \,\middle|\, \varphi \in C_c^1(\mathbb{R}^n, \mathbb{R}^n),\ \|\varphi\|_\infty \leq 1 \right\}.$$

In the special case when $\partial E$ is of class $C^1$, the perimeter $P(E)$ agrees with the natural $(n-1)$-dimensional hypersurface measure of $\partial E$. An excellent account of finite perimeter sets can be found, e.g., in [14] and [1]. The following results show that the order of density of every finite perimeter set is not less than the number
$$m_0 := n + 1^* = n + 1 + \frac{1}{n-1}$$
and more precisely that $m_0$ is the maximum order of density common to all sets of finite perimeter.

**Theorem 2.2** ([2], Lemma 4.1). *Let $E \subset \mathbb{R}^n$ be measurable and such that $P(E) < +\infty$. Then $E \sim E^{(m_0)}$.*

**Proposition 2.2** ([6], Proposition 4.1). *For all $m \in (m_0, +\infty)$ there exists a closed set $F_m \subset \mathbb{R}^n$ of positive measure and finite perimeter such that $F_m^{(m)} = \emptyset$.*

2.4. **The density-degree function.** Prior to providing the definition of the density-degree function, observe that if $E$ is a subset of $\mathbb{R}^n$ and $x \in \mathbb{R}^n$, then the set $\{k \in [n, +\infty) \,|\, x \in E^{(k)}\}$ is a (possibly empty) interval.

**Definition 2.2.** *Let $E$ be a subset of $\mathbb{R}^n$. Then define the "density-degree function of $E$" $d_E : \mathbb{R}^n \to [0, +\infty]$ as follows*
$$d_E(x) := \begin{cases} \sup \{m \in [n, +\infty) \,|\, x \in E^{(m)}\} & \text{if } x \in E^{(n)} \\ 0 & \text{if } x \notin E^{(n)}. \end{cases}$$

When there exists $k \in [n, +\infty]$ such that
$$E \sim d_E^{-1}(\{k\}) = \{x \in \mathbb{R}^n \,|\, d_E(x) = k\}$$
we say that $E$ is a "uniformly $k$-dense set".

**Example 2.1.** *If $E$ is open, then $E \subset d_E^{-1}(\{+\infty\})$. Observe that strict inclusion can occur, e.g. $E := B_r \setminus \{0\}$ (for which one has $d_E^{-1}(\{+\infty\}) = B_r$).*

**Example 2.2.** Let $m \in (2, +\infty)$ and $E$ be the set of points $(x_1, x_2) \in \mathbb{R}^2$ satisfying $|x_2| > |x_1|^{m-1}$. Since (as an elementary computation shows)
$$\lim_{r \to 0^+} \frac{\mathcal{L}^2(B(0,r) \setminus E)}{r^m} \in (0, +\infty),$$
one has $d_E(0) = m$ and $0 \notin E^{(m)}$.

This proposition collects a number of properties which have been proved in Proposition 5.1 and Proposition 5.2 of [7] (except (7) which is very easy to verify).

**Proposition 2.3.** *Let $E$ be a subset of $\mathbb{R}^n$ and $m \in [n, +\infty)$. The following properties hold:*

(1) *The density-degree function $d_E$ is measurable.*
(2) $d_E^{-1}(\{k\}) \cap d_E^{-1}(\{m\}) = \emptyset$, *if $k \neq m$ ($k \geq n$).*
(3) *If $E$ is measurable, then the set*
$$\left\{ m \in (n, +\infty) \,|\, \mathcal{L}^n(d_E^{-1}(\{m\})) > 0 \right\}$$
*is at most countable.*
(4) $d_E^{-1}((m, +\infty]) = d_E^{-1}((m, \|d_E\|_\infty]) = \bigcup_{k > m} E^{(k)}$.
(5) *If $m > n$ then*
$$d_E^{-1}([m, +\infty]) = d_E^{-1}([m, \|d_E\|_\infty]) = \bigcap_{l \in [n,m)} E^{(l)},$$
*while*
$$d_E^{-1}([n, +\infty]) = d_E^{-1}([n, \|d_E\|_\infty]) = E^{(n)}.$$
(6) $d_E^{-1}((m, +\infty]) \subset E^{(m)} \subset d_E^{-1}([m, +\infty])$.
(7) *Let $E^{(m)} \neq \emptyset$. Then $d_E|_{E^{(m)}} \equiv m$ if and only if $\|d_E\|_\infty = m$.*

**Remark 2.2.** *Both the inclusions in statement (4) of Proposition 2.3 may be strict (cf. Proposition 2.4 below and Example 2.2, respectively).*

**Remark 2.3.** *Let $\Omega$ and $E$ be, respectively, an open subset of $\mathbb{R}^n$ and a measurable subset of $\Omega$. We observe that the existence of even a single point $x \in \Omega$ such that $d_E(x) < +\infty$ yields $\mathcal{L}^n(\Omega \setminus E) > 0$. This simple observation might lead us to believe that $\mathcal{L}^n(E)$ must be small if the set of such points $x$ has a large measure. Proposition 2.4 and Theorem 2.3 below show, in particular, that this is not true.*

The following result establishes that a bounded open set in $\mathbb{R}^n$ can be arbitrarily approximated from inside by closed uniformly $n$-dense sets.

**Proposition 2.4** ([7], Proposition 5.4)**.** *Let $\Omega$ be a bounded open subset of $\mathbb{R}^n$. Then for all $C < \mathcal{L}^n(\Omega)$ there exists an uniformly n-dense closed subset $F$ of $\overline{\Omega}$ such that $\mathcal{L}^n(F) > C$.*

We expect that Proposition 2.4 can be extended to a result of approximation from inside by closed uniformly $k$-dense sets, for all $k \geq n$. We are not yet able to resolve this conjecture, but we have the following result.

**Theorem 2.3** ([7], Theorem 5.1). *Let $\Omega$ be a bounded open subset of $\mathbb{R}^n$ and let $m \in (n, +\infty)$. Then for all $C < \mathcal{L}^n(\Omega)$ and for all $t \in (n, m)$ there exist a closed subset $F$ of $\overline{\Omega}$ and an open subset $U$ of $\Omega$ such that:*

(1) $d_F^{-1}([t, +\infty]) \supset \Omega \setminus U$ and $\mathcal{L}^n(U) < \mathcal{L}^n(\Omega) - C$ *(hence $F \supset \Omega \setminus U$ and $\mathcal{L}^n(F) > C$);*
(2) *One has $F^{(m)} = \emptyset$ (hence $\|d_F\|_\infty \leq m$).*

*In particular, one has $t \leq d_F(x) \leq m$ for all $x \in \Omega \setminus U$.*

**Remark 2.4.** *When $\partial \Omega$ is Lipschitz, it is obvious that the closed set $F$ in Proposition 2.4 and in Theorem 2.3 can be chosen so that $F \subset \Omega$.*

Finally, observe that Theorem 2.2 can be restated as follows:

**Proposition 2.5** ([7], Proposition 5.3). *Let $E \subset \mathbb{R}^n$ be measurable and such that $P(E) < +\infty$. Then one has $d_E^{-1}([m_0, +\infty]) \sim E$.*

## 3. A simple characterization of superdensity points

**Proposition 3.1.** *Let $E \subset \mathbb{R}^n$ be measurable, $x \in \mathbb{R}^n$ and $m \in [n, +\infty)$.*

(1) *If $x \in E^{(m)}$, then*
$$\int_{E^c} g(y) \varphi\left(\frac{y-x}{r}\right) d\mathcal{L}^n(y) = o(r^m) \quad (as\ r \to 0+) \tag{3.1}$$
*for all $\varphi \in C_c(\mathbb{R}^n)$ and for every measurable function $g : \mathbb{R}^n \to \overline{\mathbb{R}}$ which is bounded in a neighborhood of $x$.*

(2) *Let $g \in C(\mathbb{R}^n)$ be such that $g(x) \neq 0$ and (3.1) holds for all $\varphi \in C_c(\mathbb{R}^n)$. Then $x \in E^{(m)}$.*

*Proof.* (1) Let us consider $x \in E^{(m)}$, $\varphi \in C_c(\mathbb{R}^n)$ and an arbitrary measurable function $g : \mathbb{R}^n \to \overline{\mathbb{R}}$ which is bounded in a neighborhood of $x$. If $R$ is a positive number such that $\mathrm{supp}(\varphi) \subset B(0, R)$, then

$$\left| \int_{E^c} g(y) \varphi\left(\frac{y-x}{r}\right) d\mathcal{L}^n(y) \right| = \left| \int_{B(x,rR) \cap E^c} g(y) \varphi\left(\frac{y-x}{r}\right) d\mathcal{L}^n(y) \right|$$
$$\leq \left( \sup_{B(x,rR)} |g| \right) \|\varphi\|_\infty \mathcal{L}^n(B(x, rR) \cap E^c).$$

Hence (3.1) follows at once.

(2) Let us consider $\varphi \in C_c(\mathbb{R}^n)$ such that
$$\varphi(\mathbb{R}^n) = [0,1], \quad \mathrm{supp}(\varphi) \subset B(0,2), \quad \varphi|_{B(0,1)} = 1.$$
Moreover, without loss of generality, we can assume that $g$ is positive in a neighborhood of $x$. Hence two positive constants $p$ and $r_0$ have to exist such that $\inf_{B(x,2r_0)} g \geq p$. For all $r \in (0, r_0]$ we have
$$p\,\mathcal{L}^n(B(x,r) \cap E^c) \leq \int_{B(x,r) \cap E^c} g(y)\,d\mathcal{L}^n(y)$$
$$\leq \int_{E^c} g(y)\,\varphi\left(\frac{y-x}{r}\right) d\mathcal{L}^n(y),$$
hence the conclusion follows from (3.1). □

This corollary is a trivial consequence of Proposition 3.1.

**Corollary 3.1.** *Let $E \subset \mathbb{R}^n$ be measurable, $x \in \mathbb{R}^n$ and $m \in [n, +\infty)$. Then each of the following properties is equivalent to $x \in E^{(m)}$:*

(1) *The equation*
$$\int_{\frac{E-x}{r}} \varphi\,d\mathcal{L}^n = \int_{\mathbb{R}^n} \varphi\,d\mathcal{L}^n + o(r^{m-n}) \qquad (\text{as } r \to 0+)$$
*holds for all $\varphi \in C_c(\mathbb{R}^n)$.*
(2) *There exists $g \in C(\mathbb{R}^n)$ such that $g(x) \neq 0$ and (3.1) holds for all $\varphi \in C_c(\mathbb{R}^n)$.*
(3) *The identity (3.1) holds for all $\varphi \in C_c(\mathbb{R}^n)$ and for every measurable function $g : \mathbb{R}^n \to \overline{\mathbb{R}}$ which is bounded in a neighborhood of $x$.*

The next result is also a very easy consequence of Proposition 3.1.

**Corollary 3.2.** *Let us consider a measurable set $E \subset \mathbb{R}^n$, $x \in E^{(m)}$ with $m \in [n, +\infty)$ and a measurable function $\Gamma : \mathbb{R}^n \to \overline{\mathbb{R}}$ which is continuous at $x$. Assume that there exists $\psi \in C_c(\mathbb{R}^n)$ such that $\int_{\mathbb{R}^n} \psi\,d\mathcal{L}^n \neq 0$ and (for $r$ small enough)*
$$\left|\int_{\mathbb{R}^n} \Gamma(y)\psi\left(\frac{y-x}{r}\right) d\mathcal{L}^n(y)\right| \leq r^{n-m} \sum_{i=1}^k \left|\int_{E^c} g_i(y)\varphi_i\left(\frac{y-x}{r}\right) d\mathcal{L}^n(y)\right|$$
*where $g_1, \ldots, g_k : \mathbb{R}^n \to \overline{\mathbb{R}}$ is a family of measurable functions which are bounded in a neighborhood of $x$ and $\varphi_1, \ldots, \varphi_k \in C_c(\mathbb{R}^n)$. Then $\Gamma(x) = 0$.*

*Proof.* From (3.1) it follows that
$$\int_{\mathbb{R}^n} \Gamma(y)\psi\left(\frac{y-x}{r}\right) d\mathcal{L}^n(y) = r^{n-m} o(r^m) = o(r^n) \qquad (\text{as } r \to 0+).$$

On the other hand, we have also

$$\int_{\mathbb{R}^n} \Gamma(y)\psi\left(\frac{y-x}{r}\right) d\mathcal{L}^n(y) = r^n \int_{\mathbb{R}^n} \Gamma(x+rz)\psi(z) \, d\mathcal{L}^n(z),$$

so that

$$\int_{\mathbb{R}^n} \Gamma(x+rz)\psi(z) \, d\mathcal{L}^n(z) = r^{-n} o(r^n) \qquad (\text{as } r \to 0+).$$

Hence (letting $r \to 0+$)

$$\Gamma(x) \int_{\mathbb{R}^n} \psi \, d\mathcal{L}^n = 0$$

and the conclusion follows recalling that $\int_{\mathbb{R}^n} \psi \, d\mathcal{L}^n \neq 0$. □

We are interested in Corollary 3.2 as it provides a common argument for the proofs of several theorems we have obtained in our work on superdensity, e.g., [2, Theorem 2.1] (the oldest one) and [10, Theorem 3.6] (the most recent one). Having it in mind can be useful for yielding new applications. As an example, let us prove the following result.

**Theorem 3.1.** *Let $\lambda$ be a continuous differential form of degree $h$ on $\mathbb{R}^n$ (with $h \geq 0$) having the following property: There exists a continuous differential form $\Delta$ of degree $h+1$ in $\mathbb{R}^n$ such that*

(3.2) $$\int_{\mathbb{R}^n} \Delta \wedge \omega = \int_{\mathbb{R}^n} \lambda \wedge d\omega,$$

*for every $C_c^\infty$ differential form $\omega$ of degree $n - h - 1$ in $\mathbb{R}^n$. Moreover let $\mu$ be a $C^1$ differential form of degree $h + 1$ in $\mathbb{R}^n$ and set $E := \{y \in \mathbb{R}^n \,|\, \Delta(y) = \mu(y)\}$. Then $(d\mu)|_{E^{(n+1)}} \equiv 0$.*

*Proof.* Let $x \in E^{(n+1)}$. Consider an arbitrary $C^\infty$ differential form $\theta$ of degree $n - h - 2$ in $\mathbb{R}^n$ and define $\Gamma : \mathbb{R}^n \to \mathbb{R}$ as

$$\Gamma(y) := \langle d\mu(y) \wedge \theta(y), dx_1 \wedge \cdots \wedge dx_n \rangle \qquad (y \in \mathbb{R}^n).$$

Moreover let $\psi \in C_c^\infty(\mathbb{R}^n)$ be such that $\int_{\mathbb{R}^n} \psi d\mathcal{L}^n \neq 0$. Then, from the Stokes theorem, we obtain (for $r$ small enough)

$$\begin{aligned}
\int_{\mathbb{R}^n} \Gamma(y)\psi\left(\frac{y-x}{r}\right) d\mathcal{L}^n(y) &= \int_{\mathbb{R}^n} \psi\left(\frac{y-x}{r}\right) d\mu(y) \wedge \theta(y) \\
&= -\int_{\mathbb{R}^n} d\left[\psi\left(\frac{y-x}{r}\right)\right] \wedge \mu(y) \wedge \theta(y) \\
&\quad + (-1)^h \int_{\mathbb{R}^n} \psi\left(\frac{y-x}{r}\right) \mu(y) \wedge d\theta(y) \\
&= -\int_E d\left[\psi\left(\frac{y-x}{r}\right)\right] \wedge \Delta(y) \wedge \theta(y) \\
&\quad - \int_{E^c} d\left[\psi\left(\frac{y-x}{r}\right)\right] \wedge \mu(y) \wedge \theta(y) \\
&\quad + (-1)^h \int_E \psi\left(\frac{y-x}{r}\right) \Delta(y) \wedge d\theta(y) \\
&\quad + (-1)^h \int_{E^c} \psi\left(\frac{y-x}{r}\right) \mu(y) \wedge d\theta(y) \\
&= -\int_{\mathbb{R}^n} d\left[\psi\left(\frac{y-x}{r}\right)\right] \wedge \Delta(y) \wedge \theta(y) \\
&\quad + \int_{E^c} d\left[\psi\left(\frac{y-x}{r}\right)\right] \wedge [\Delta(y) - \mu(y)] \wedge \theta(y) \\
&\quad + (-1)^h \int_{\mathbb{R}^n} \psi\left(\frac{y-x}{r}\right) \Delta(y) \wedge d\theta(y) \\
&\quad + (-1)^h \int_{E^c} \psi\left(\frac{y-x}{r}\right) [\mu(y) - \Delta(y)] \wedge d\theta(y).
\end{aligned}$$

Hence

(3.3) $$\int_{\mathbb{R}^n} \Gamma(y)\psi\left(\frac{y-x}{r}\right) d\mathcal{L}^n(y) = I(r) + J(r),$$

where

$$I(r) := \int_{\mathbb{R}^n} -d\left[\psi\left(\frac{y-x}{r}\right)\right] \wedge \Delta(y) \wedge \theta(y) + (-1)^h \psi\left(\frac{y-x}{r}\right) \Delta(y) \wedge d\theta(y)$$

and

$$J(r) := \int_{E^c} d\left[\psi\left(\frac{y-x}{r}\right)\right] \wedge [\Delta(y) - \mu(y)] \wedge \theta(y) + (-1)^h \psi\left(\frac{y-x}{r}\right) [\mu(y) - \Delta(y)] \wedge d\theta(y).$$

Now observe that

$$\begin{aligned}
(-1)^h I(r) &= \int_{\mathbb{R}^n} \Delta(y) \wedge \left(d\left[\psi\left(\frac{y-x}{r}\right)\right] \wedge \theta(y) + \psi\left(\frac{y-x}{r}\right) d\theta(y)\right) \\
&= \int_{\mathbb{R}^n} \Delta(y) \wedge d\left[\psi\left(\frac{y-x}{r}\right) \theta(y)\right],
\end{aligned}$$

hence

(3.4) $$I(r) = 0,$$

by assumption (3.2). Moreover, for all $r \in (0, 1]$, one has

$$|J(r)| = \left| \frac{1}{r} \sum_{i=1}^{n} \int_{E^c} (D_i \psi) \left( \frac{y-x}{r} \right) dx_i \wedge [\Delta(y) - \mu(y)] \wedge \theta(y) \right.$$

(3.5)
$$\left. + (-1)^h \int_{E^c} \psi \left( \frac{y-x}{r} \right) [\mu(y) - \Delta(y)] \wedge d\theta(y) \right|$$

$$\leq \frac{1}{r} \sum_{i=1}^{n+1} \left| \int_{E^c} g_i(y) \varphi_i \left( \frac{y-x}{r} \right) d\mathcal{L}^n(y) \right|$$

with

$$\varphi_i := D_i \psi, \quad g_i := \langle dx_i \wedge (\Delta - \mu) \wedge \theta, dx_1 \wedge \cdots \wedge dx_n \rangle \quad (i = 1, \ldots, n)$$

and

$$\varphi_{n+1} := \psi, \quad g_{n+1} := \langle (\mu - \Delta) \wedge d\theta, dx_1 \wedge \cdots \wedge dx_n \rangle.$$

From (3.3), (3.4), (3.5) and Corollary 3.2 (with $m = k = n + 1$) we obtain $\Gamma(x) = 0$. The conclusion follows from the arbitrariness of $\theta$, by recalling (2.2). □

**Remark 3.1.** Let $\lambda$ and $\omega$ be, respectively, a $C^1$ differential form of degree $h$ and a $C_c^1$ differential form of degree $n - h - 1$ on $\mathbb{R}^n$. Since $d(\lambda \wedge \omega) = d\lambda \wedge \omega + (-1)^h \lambda \wedge d\omega$, we have

$$\int_{\mathbb{R}^n} \lambda \wedge d\omega = (-1)^h \int_{\mathbb{R}^n} d(\lambda \wedge \omega) + (-1)^{h+1} \int_{\mathbb{R}^n} d\lambda \wedge \omega = (-1)^{h+1} \int_{\mathbb{R}^n} d\lambda \wedge \omega.$$

So, when $\lambda$ is of class $C^1$, the condition (3.2) is verified with $\Delta = (-1)^{h+1} d\lambda$. Hence Theorem 3.1 yields the following result (provided in [2, Theorem 2.1]): Let $\lambda$ and $\mu$ be $C^1$ differential forms on $\mathbb{R}^n$ of degree $h$ and $h+1$, respectively (with $h \geq 0$). If define $E := \{y \in \mathbb{R}^n \,|\, d\lambda(y) = \mu(y)\}$, then $(d\mu)|_{E^{(n+1)}} \equiv 0$.

## 4. Some further properties of the density-degree functions

**Proposition 4.1.** Given two subsets $E, F$ of $\mathbb{R}^n$, the following properties hold:

(1) If $E \sim F$, then $E^{(m)} = F^{(m)}$ for all $m \in [n, +\infty)$, hence $d_E = d_F$.
(2) If $E, F$ are measurable and $E \not\sim F$, then $\{d_E \neq d_F\}$ has positive measure.

*Proof.* (1) Since $\mathcal{L}^n(E \cap F^c) = \mathcal{L}^n(F \cap E^c) = 0$, the sets $E \cap F^c$ and $F \cap E^c$ are measurable. It follows that (for all $x \in \mathbb{R}^n$ and $r > 0$)

$$\mathcal{L}^n(B(x,r) \setminus E) = \mathcal{L}^n(B(x,r) \cap [(E \cap F) \cup (E \cap F^c)]^c)$$
$$= \mathcal{L}^n(B(x,r) \cap (E \cap F)^c \cap (E \cap F^c)^c)$$
$$= \mathcal{L}^n(B(x,r) \cap (E \cap F)^c) - \mathcal{L}^n(B(x,r) \cap (E \cap F)^c \cap (E \cap F^c))$$

where the second term on the right above is 0 by monotonicity. Since the last term is symmetric in $E, F$, one has

$$\mathcal{L}^n(B(x,r) \setminus F) = \mathcal{L}^n(B(x,r) \setminus E) = \mathcal{L}^n(B(x,r) \cap (E \cap F)^c),$$

hence $E^{(m)} = F^{(m)}$ for all $m \geq n$.

(2) First of all, recall that (by a well-known result on the Lebesgue set, e.g., see Corollary 1.5 in [18, Chapter 3]) a measure zero set $N \subset \mathbb{R}^n$ has to exist such that

(4.1) $$E \setminus N \subset E^{(n)}, \quad (F^c \setminus N) \cap F^{(n)} = \emptyset$$

and

(4.2) $$F \setminus N \subset F^{(n)}, \quad (E^c \setminus N) \cap E^{(n)} = \emptyset.$$

Moreover, by hypothesis, at least one of the following inequalities must hold:

$$\mathcal{L}^n(E \cap F^c) > 0, \quad \mathcal{L}^n(F \cap E^c) > 0.$$

- If the first inequality holds, then, by (4.1), for all $x$ in the positive measure set $E \cap F^c \cap N^c$ one has $d_E(x) \geq n$ (in that $x \in E \cap N^c \subset E^{(n)}$) and $d_F(x) = 0$ (in that $x \in F^c \cap N^c$, hence $x \notin F^{(n)}$).
- If instead the second inequality holds, then, by (4.2), for all $x$ in the positive measure set $F \cap E^c \cap N^c$ one has $d_F(x) \geq n$ (in that $x \in F \cap N^c \subset F^{(n)}$) and $d_E(x) = 0$ (in that $x \in E^c \cap N^c$, hence $x \notin E^{(n)}$).

$\square$

The following property follows immediately from (1) of Proposition 4.1, by also recalling that Lebesgue outer measure is Borel regular.

**Corollary 4.1.** *For every measurable subset $E$ of $\mathbb{R}^n$ there exists a Borel set $B$ such that $E \subset B$, $\mathcal{L}^n(E) = \mathcal{L}^n(B)$ and $d_E = d_B$.*

**Proposition 4.2.** *Let $\{E_j\}_{j \in J}$ be any family of subsets of $\mathbb{R}^n$. Then the following inequalities hold:*

(4.3) $$d_{\cap_{j \in J} E_j} \leq \inf_{j \in J} d_{E_j}, \quad d_{\cup_{j \in J} E_j} \geq \sup_{j \in J} d_{E_j}.$$

*If $J$ is finite, then the first one turns into the equality $d_{\cap_{j \in J} E_j} = \min_{j \in J} d_{E_j}$, while the identity $d_{\cup_{j \in J} E_j} = \max_{j \in J} d_{E_j}$ fails to be true in general.*

*Proof.* For all $x \notin (\cap_{j \in J} E_j)^{(n)}$ we obviously have

$$d_{\cap_{j \in J} E_j}(x) = 0 \leq \inf_{j \in J} d_{E_j}(x).$$

So let us suppose $x \in (\cap_{j \in J} E_j)^{(n)}$. Then the set $\{k \geq n \mid x \in (\cap_{j \in J} E_j)^{(k)}\}$ is non-empty and
$$\begin{aligned} d_{\cap_{j \in J} E_j}(x) &= \sup\{k \geq n \mid x \in (\cap_{j \in J} E_j)^{(k)}\} \\ &\leq \sup\{k \geq n \mid x \in \cap_{j \in J} E_j^{(k)}\} \\ &= \inf_{j \in J} \sup\{k \geq n \mid x \in E_j^{(k)}\} \\ &= \inf_{j \in J} d_{E_j}(x), \end{aligned}$$
where equality holds in the second line whenever $J$ is finite (cf. Remark 2.1). This concludes the proof of the part concerning $d_{\cap_{j \in J} E_j}$.

For all $x \notin (\cup_{j \in J} E_j)^{(n)}$, one also has $x \notin E_j^{(n)}$ for all $j \in J$ (cf. Remark 2.1) and thus
$$d_{\cup_{j \in J} E_j}(x) = \sup_{j \in J} d_{E_j}(x) = 0.$$

On the other hand, if $x \in (\cup_{j \in J} E_j)^{(n)}$, then the set $\{k \geq n \mid x \in (\cup_{j \in J} E_j)^{(k)}\}$ is non-empty and (cf. Remark 2.1)
$$\begin{aligned} d_{\cup_{j \in J} E_j}(x) &= \sup\{k \geq n \mid x \in (\cup_{j \in J} E_j)^{(k)}\} \\ &\geq \sup\{k \geq n \mid x \in \cup_{j \in J} E_j^{(k)}\} \\ &= \sup_{j \in J} \sup\{k \geq n \mid x \in E_j^{(k)}\} \\ &= \sup_{j \in J} d_{E_j}(x). \end{aligned}$$

To verify the last assertion, we can consider the example provided in Remark 2.1 ($n = 1$, $J = \{1, 2\}$, $E_1 = (-1, 0)$, $E_2 = (0, 1)$) which yields
$$d_{\cup_{j \in J} E_j} = +\infty \chi_{(-1,1)}, \quad \max_{j \in J} d_{E_j} = +\infty \chi_{(-1,1) \setminus \{0\}}.$$
□

## 5. Does the set having an arbitrarily prescribed density-degree function always exist?

As we have stated in Proposition 2.3, every density-degree function (relative to a subset of $\mathbb{R}^n$) is measurable and takes its values in $\{0\} \cup [n, +\infty]$. Rough considerations on this subject could lead us to believe that for every measurable function $f : \mathbb{R}^n \to \overline{\mathbb{R}}$ such that $f(x) \in \{0\} \cup [n, +\infty]$ at a.e. $x \in \mathbb{R}^n$ there exists $E \subset \mathbb{R}^n$ for which $d_E \sim f$. The following example disproves this.

**Example 5.1.** *Let $F$ be a measurable subset of $\mathbb{R}^n$ such that $0 < \|d_F\|_\infty < +\infty$ (observe that $F$ can be provided by Proposition 2.4 or Theorem 2.3) and let $m \in (\|d_F\|_\infty, +\infty)$. Then consider the measurable function $f := m\chi_F$ and prove (by contradiction) that there is no measurable set $E \subset \mathbb{R}^n$ for which $d_E \sim f$. In fact, if such a set $E$ exists, then* $E \sim E^{(n)} \sim F$ (cf. Definition 2.2), hence $d_F = d_E \sim m\chi_F$ (cf. Proposition 4.1) and this contradicts the hypothesis $m > \|d_F\|_\infty$.

These limitations, however, are not sufficient to prevent good approximation properties. In fact, the following result holds.

**Theorem 5.1.** *Let $f : \mathbb{R}^n \to \overline{\mathbb{R}}$ be a measurable function such that $f(x) \in \{0\} \cup [n, +\infty]$ for a.e. $x \in \mathbb{R}^n$. Then there exists a countable family $\{F_k\}_{k=1}^\infty$ of closed subsets of $\mathbb{R}^n$ such that $\lim_{k \to \infty} d_{F_k}(x) = f(x)$ for a.e. $x \in \mathbb{R}^n$.*

*Proof.* Let $N := f^{-1}([-\infty, 0) \cup (0, n))$ and define
$$\widetilde{f}(x) := \begin{cases} f(x) & \text{if } x \in \mathbb{R}^n \setminus N \\ 0 & \text{if } x \in N. \end{cases}$$

Note that $\widetilde{f}$ is equivalent to $f$ (since $\mathcal{L}(N) = 0$, by assumption) and takes all of its values in $\{0\} \cup [n, +\infty]$. By [17, Theorem 1.17] there exists a nondecreasing sequence of simple measurable functions on $\mathbb{R}^n$
$$s_k = \sum_{j=1}^{M_k} a_{k,j} \chi_{E_{k,j}} \qquad (k = 1, 2, \dots)$$
with the following properties:

(i) $\{E_{k,j}\}_{j=1}^{M_k}$ is a measurable partition of $E := \widetilde{f}^{-1}([n, +\infty])$ and $\{a_{k,j}\}_{j=1}^{M_k} \subset [n, +\infty)$;
(ii) $\lim_{k \to \infty} s_k(x) = \widetilde{f}(x)$, for all $x \in \mathbb{R}^n$.

Now arguing as in the proof of [18, Theorem 4.3], with
$$\varphi_k := s_k \chi_{B_k} = \sum_{j=1}^{M_k} a_{k,j} \chi_{E_{k,j} \cap B_k} \qquad (k = 1, 2, \dots),$$
we can find a sequence of simple functions
$$\psi_k = \sum_{j=1}^{N_k} b_{k,j} \chi_{R_{k,j}} \qquad (k = 1, 2, \dots),$$
where $\{R_{k,j}\}_{j=1}^{N_k}$ is a collection of disjoint closed rectangles and $\{b_{k,j}\}_{j=1}^{N_k} \subset [n, +\infty)$, such that
(5.1) $$\lim_{k \to \infty} \psi_k(x) = \widetilde{f}(x) = f(x), \text{ for a.e. } x \in \mathbb{R}^n.$$

Without loss of generality we can also assume $b_{k,j} \in (n, +\infty)$ for all $k, j$ (it is enough to replace $b_{k,j}$ with $b_{k,j} + 1/k$), hence a sequence of positive real numbers $\{\varepsilon_k\}_{k=1}^\infty$ has to exist such that $\lim_{k \to \infty} \varepsilon_k = 0$ and
$$b_{k,j} - \varepsilon_k > n, \text{ for all } j = 1, \dots, N_k.$$

Now let $\Omega_{k,j}$ denote the interior of $R_{k,j}$. By applying Theorem 2.3 and Remark 2.4, we find a closed set $F_{k,j} \subset \Omega_{k,j}$ and an open set $U_{k,j} \subset \Omega_{k,j}$ satisfying

$$(5.2) \qquad \mathcal{L}^n(U_{k,j}) < \frac{2^{-k}}{N_k}, \quad \mathcal{L}^n(F_{k,j}) > \mathcal{L}^n(\Omega_{k,j}) - \frac{2^{-k}}{N_k}$$

and

$$|d_{F_{k,j}}(x) - b_{k,j}| < \varepsilon_k, \text{ for all } x \in \Omega_{k,j} \setminus U_{k,j}.$$

In particular, if define

$$(5.3) \qquad \Omega_k := \bigcup_{j=1}^{N_k} \Omega_{k,j}, \quad F_k := \bigcup_{j=1}^{N_k} F_{k,j}, \quad U_k := \bigcup_{j=1}^{N_k} U_{k,j},$$

then we get

$$(5.4) \qquad |d_{F_k} - \psi_k| < \varepsilon_k \text{ in } \Omega_k \setminus U_k.$$

Also one obviously has

$$(5.5) \qquad d_{F_k} = \psi_k = 0 \text{ in } R_k^c, \text{ where } R_k := \bigcup_{j=1}^{N_k} R_{k,j}.$$

Observe that

$$W := \left(\bigcap_{l=1}^{\infty} \bigcup_{k>l} U_k\right) \cup \left(\bigcup_{h=1}^{\infty} \partial R_h\right)$$

is a set of measure zero, by (5.2) and (5.3). Moreover, if

$$x \in W^c = \left(\bigcup_{l=1}^{\infty} \bigcap_{k>l} U_k^c\right) \cap \left(\bigcup_{h=1}^{\infty} \partial R_h\right)^c = \bigcup_{l=1}^{\infty} \left[\left(\bigcap_{k>l} U_k^c\right) \cap \left(\bigcup_{h=1}^{\infty} \partial R_h\right)^c\right],$$

then there exists $l_x \geq 1$ such that

$$x \in \left(\bigcap_{k>l_x} U_k^c\right) \cap \left(\bigcup_{h=1}^{\infty} \partial R_h\right)^c,$$

namely

$$x \in U_k^c \cap \left(\bigcup_{h=1}^{\infty} \partial R_h\right)^c \subset (\Omega_k \setminus U_k) \cup R_k^c, \text{ for all } k > l_x.$$

Hence, by (5.4) and (5.5), we get

$$|d_{F_k}(x) - \psi_k(x)| < \varepsilon_k \text{ for all } k > l_x.$$

The conclusion follows from (5.1). □

**Corollary 5.1.** *Let $E$ be a measurable subset of $\mathbb{R}^n$. Then there exists an increasing sequence $\{F_l\}_{l=1}^{\infty}$ of uniformly $n$-dense closed subsets of $\mathbb{R}^n$ such that $\cup_l F_l \sim E$.*

*Proof.* By applying Theorem 5.1 to the measurable function $f := n\chi_E$, we find a countable family $\{\widetilde{F}_k\}_{k=1}^{\infty}$ of closed subsets of $\mathbb{R}^n$ and a measure zero set $N \subset \mathbb{R}^n$ such that

$$\lim_{k\to\infty} d_{\widetilde{F}_k}(x) = n, \text{ for all } x \in E \setminus N \tag{5.6}$$

and

$$\lim_{k\to\infty} d_{\widetilde{F}_k}(x) = 0, \text{ for all } x \in E^c \setminus N. \tag{5.7}$$

From (5.6) and (5.7) we get, respectively (cf. Definition 2.2),

$$E \setminus N \subset \bigcup_l \bigcap_{k\geq l} \widetilde{F}_k^{(n)}, \qquad E^c \setminus N \subset \bigcup_l \bigcap_{k\geq l} (\widetilde{F}_k^{(n)})^c. \tag{5.8}$$

From the second inclusion in (5.8) it follows that

$$\bigcup_l \bigcap_{k\geq l} \widetilde{F}_k^{(n)} \subset \bigcap_l \bigcup_{k\geq l} \widetilde{F}_k^{(n)} \widetilde{\subset} E.$$

Hence and by the first inclusion in (5.8) we obtain

$$E \sim \bigcup_l \bigcap_{k\geq l} \widetilde{F}_k^{(n)} \sim \bigcup_l F_l, \text{ with } F_l := \bigcap_{k\geq l} \widetilde{F}_k.$$

Observe that $\{F_l\}_{l=1}^{\infty}$ is an increasing sequence of closed sets. The only thing left to prove is that each $F_l$ is uniformly $n$-dense. For this purpose, observe first that for all $x \in F_l^{(n)} \cap E \setminus N$ we have

$$n \leq d_{F_l}(x) \leq \inf_{k\geq l} d_{\widetilde{F}_k}(x) = n,$$

by Definition 2.2, Proposition 4.2 and (5.6). Now the conclusion follows from the equivalence $F_l^{(n)} \cap E \setminus N \sim F_l$. □